\documentclass{amsart}
\usepackage{amsmath,amssymb}
\usepackage[dvips]{epsfig}
\usepackage[all, knot]{xy}
\usepackage{hyperref}
\SelectTips{cm}{}
\xyoption{arc}
\xyoption{2cell}
\UseAllTwocells
\newcommand \driicell [1] {\drtwocell \omit {#1}}

\newcommand \drriicell [1] {\drrtwocell \omit {#1}}
\newcommand \uriicell [1] {\urtwocell \omit {#1}}

\newtheorem{theorem}{Theorem}
\numberwithin{theorem}{subsection}

\newtheorem{varexample}[theorem]{Example}

\newtheorem{varremark}[theorem]{Remark}
\newenvironment{remark}{\begin{varremark}\em}{\em\end{varremark}}

\newcommand{\opname}[1]{\operatorname{{#1}}}
\newcommand{\catname}[1]{\boldsymbol{\operatorname{{#1}}}}
\newcommand{\C}{\catname{C}}
\newcommand{\Cat}{\catname{Cat}}
\newcommand{\Set}{\catname{Set}}
\newcommand{\Bicat}{\catname{Bicat}}
\newcommand{\Hor}{\catname{Hor}}
\newcommand{\Ver}{\catname{Ver}}

\newcommand{\ra}{\rightarrow}

\newtheorem{definition}{Definition}

\begin{document}
\title{Cubical $n$-Categories and Finite-Limits Theories}

\author{Jeffrey C. Morton}

\begin{abstract}
This note informally describes a way to build certain cubical
$n$-categories by iterating a process of taking models of certain
finite limits theories.  We base this discussion on a construction of
``double bicategories'' as bicategories internal to $\Bicat$, and see
how to extend this to $n$-tuple bicategories (and similarly for
tricategories etc.)  We briefly consider how to reproduce ``simpler''
definitions of weak cubical $n$-category from these.
\end{abstract}

\maketitle

\section{Introduction}

This note aims to describe, through a particular example, some
relationships between two of the main topics of MakkaiFest: model
theory and higher categories. In particular, we aim to describe a way
to build certain higher categories by iterating a process of taking
models of a theory. The higher categories we have in mind are a
particular sort of weak cubical $n$-category.

While we make no attempt to be comprehensive, we will start with a
particular example, an extension of the view of double categories as
categories internal to $\Cat$.  Double categories are structures which
have two different category structures on the same set of objects, one
``horizontal'' and one ``vertical''.  Often these directions are
different, such as the double category whose objects are sets, whose
horizontal morphisms are functions, and whose vertical morphisms are
relations.  There are also ``squares'' which act as non-identity cells
which can fill square diagrams:
\begin{equation}
\xymatrix{
X \driicell{F} \ar[r]^{g} \ar[d]_{\xi} & Y \ar[d]^{\psi} \\
X' \ar[r]_{g'} & Y'
}
\end{equation}

We take as starting point ``double bicategories'', found by taking
bicategories internal to $\Bicat$, where composition satisfies weaker
axioms; in particular there is an associator isomorphism for each
composition.  This structure turns out to be useful for describing
structures build from cospans; in particular, a weak cubical
2-category of cobordisms between cobordisms between manifolds, in any
dimension.  These are most naturally taken with weak composition,
since composition is by gluing of smooth manifolds (by a
diffeomorphism).

The author has described the resulting structure elsewhere
\cite{dblbicat}, and related it to a related definition given by
Dominic Verity.  Verity's double bicategories can be obtained from
internal bicategories in $\Bicat$ in cases where a certain
niche-filling condition obtains.  We will discuss here how this can be
specified as a property of a model, and thus find more restricted
notions of cubical $n$-category as special cases.

First we recap the idea of finite limits theories, and how categories
and bicategories, at least (and similarly globular $n$-categories for
any $n$) are described in terms of finite limits theories.  We note
how to obtain double categories, or double bicategories respectively,
and describe how this construction generalizes by iterated application
of the process of taking models.

That is, given a finite limits theory, say $\catname{Th(Cat)}$, the
\textit{theory of categories}, one repeatedly applies the functor
$[\catname{Th(Cat)}, -]$.  This notation indicates the category of
models of $\catname{Th(Cat)}$ in the target category.  Thus, one takes
categories of models, and then finds models in \textit{these}
categories.  The structures formed in this was are strict cubical
$n$-categories.  Doing the same with the theory of bicategories,
$\catname{Th(Bicat)}$, gives the particular notion of weak cubical
$n$-category we are discussing, namely, $n$-tuple bicategories.  We
will not discuss tricategories, tetracategories, or other forms of
$n$-category here, but note a similar process gives $n$-tuple
tricategories etc.

There are other possible generalizations of this way of generating
cubical $n$-categories.  One involves weakening the notion of model:
we will not consider this in depth, but will consider the example of
pseudocategories to see how this approach motivated the one we are
taking here.  We discuss will double categories and pseudocategories
as strict and weak models of the theory $\catname{Th(Cat)}$, and how
this motivates double bicategories as models of
$\catname{Th(\Bicat)}$.  We will not consider, here, the full extent
of weakening possible due to the fact that $\Bicat$ is not just a
category, but a tricategory.

We also briefly consider the ``niche-filling'' conditions which can be
used to produce Verity double bicategories, which have fewer different
types of morphism.  It is possible to extend this to higher $n$, but
the conditions become more complex, and we will regard it as more
elegant to leave the extra types of morphism in place, and conclude by
indicating an example where this is quite natural.

\section{The Theory of Bicategories}\label{sec:thcatbicat}

To begin with, we recall that, following the approach of Lawvere
\cite{wl} to universal algebra (and see also \cite{borceux2}), an
algebraic theory can be understood as a category. For example, the
theory of groups can be seen as $\mathbf{Th(Grp)}$, the free Cartesian
category on a group object. That is, $\mathbf{Th(Grp)}$ is the minimal
category with finite products, a terminal object, and an object $G$
equipped with maps $m : G \times G \rightarrow G$ and $e : 1
\rightarrow G$ satisfying the usual group axioms. For example, $(m
\otimes 1)\circ m = m \circ (1 \otimes m)$.  Then a product-preserving
functor $M : \mathbf{Th(Grp)} \rightarrow \mathcal{C}$, into a
Cartesian category $\mathcal{C}$ amounts to the same thing as a
``group object in $\mathcal{C}$'' in the usual sense. That is, the
group object is $M(G)$, and the multiplication and unit maps are
$M(m)$ and $M(e)$. In the case that $\mathcal{C} = \mathbf{Sets}$,
such a functor is just a group.

A theory for a structure $\mathcal{S}$ is described by a category
$\catname{Th}(\mathcal{S})$ which we think of as a \textit{diagram}
containing all the axioms defining $\mathcal{S}$.  A \textit{model} of
such a theory in a category $\catname{C}$ is a functor into $\catname{C}$.

Some structures cannot be described by \textit{algebraic} theories,
however. These categories need only have objects $\{ T^0, T^1, T^2,
\dots \}$, the powers of a given object, such as the object $G \in
\catname{Th(Grp)}$.  Taking models $F : \catname{Th(Grp)} \ra
\catname{Set}$, one then gets the usual definition of a group as the
set $F(G)$, equipped with set maps such as $F(m) : F(G) \times F(G)
\ra F(G)$, satisfying the axioms.  However, some
structures---categories themselves, for example---are not naturally
described in terms of a single set (i.e. a single object in the
theory, and its powers). This raises the question of the
\textit{doctrine} of the theory---that is, the specific 2-category of
categories in which we take the theory, and the functors which are its
models. The doctrine specifies what structures are required to define
the theory (as products, or at least monoidal structure, are required
to define the multiplication map for the theory of groups).  These are
the structures which must be preserved by the model functors.

To begin describing $n$-categories in terms of models, we need to
describe the \textit{theory of categories}, and indeed of
bicategories.  The crucial difference between the (algebraic) theory
of groups, and the theory of categories is that the composition
operation for morphisms, which plays the role of the multiplication
map $m$ in the group example, is only partially defined.  This means
$\catname{Th(Cat)}$ must have some extra structure.

\begin{definition}
The doctrine of categories with finite limits is a 2-category
$\mathbf{FL}$ whose objects are categories with all finite limits, as
morphisms all functors which preserve finite limits, and as
2-morphisms all natural transformations.  A finite limits theory is a
category $\mathcal{T}$ in $\mathbf{FL}$, and a model of $\mathcal{T}$
is a functor $\mathcal{T} \rightarrow \mathcal{C}$.
\end{definition}

(Of course, we only usually emphasize that a category is a theory if
it is easily presentable in terms of generators and relations, as in
the theory of groups.) 

We want to describe a whole class of concepts of cubical $n$-category
which can be framed in terms of model theory.  These are (some sort
of) models of a theory of (globular, or other) $n$-categories in a
suitable category of such $n$-categories.  There are many variants,
depending on the choice of theory and the target category for model.
We want a ``weak'' notion, in the sense that composition is
associative up to an isomorphism, rather than exactly.  The most
elementary higher categories with this property are bicategories.

To begin with, we consider the simpler theory of categories,
$\catname{Th(Cat)}$.  The usual definition of a category will amount
to a model of $\catname{Th(Cat)}$ in $\Set$.

\begin{definition}The \textit{theory} $\catname{Th(Cat)}$ is a
category in $\textbf{FL}$, generated by the following data:
\begin{itemize} 
\item two objects $Obj$ and $Mor$ 
\item morphisms of the form:
\begin{equation}
  \xymatrix{
    Mor \ar@/^/[r]^{s} \ar@/_/[r]_{t} & Obj \\
  }
\end{equation}
and
\begin{equation}
  \xymatrix{
    Obj \ar[r]^{id} & Mor  \\
  }
\end{equation} such that $s (id) = t (id) = 1_{Obj}$ as expected.
\item if $\opname{Pairs}$ is the pullback in the square:
\begin{equation}\label{eq:pairs}
  \xymatrix{
     & \opname{Pairs} \ar[dl]^{p_1} \ar[dr]_{p_2} & \\
     Mor \ar[dr]^{t} &  & Mor \ar[dl]_{s} & \\
     & Obj  & \\
  }
\end{equation} then there is a (``partially defined'') composition map
\[
\circ : \opname{Pairs} \rightarrow Mor
\]
satisfying the usual properties for composition, namely:
\begin{equation}\label{eq:pullbackcomposn}
  \xymatrix{
      &   & Mor \ar@/_2pc/[dddll]_{s} \ar@/^2pc/[dddrr]^{t} & & \\
      &   & \opname{Pairs} \ar[u]^{\circ} \ar[dl]_{p_1} \ar[dr]^{p_2}  &  & \\
      & Mor \ar[dl]_{s} \ar[dr]^{t} &   & Mor \ar[dl]_{s} \ar[dr]^{t} & \\
    Obj &  & Obj &   & Obj \\
  }
\end{equation} (and another axiom with the interpretation in
$\catname{Sets}$ that for any morphism $f \in Mor$, we have
$id(s(f))$ and $1(t(f))$ are composable with $f$, and the composite is
$f$).
\end{itemize}
\end{definition}

\begin{remark}For categories in $\Set$, the pullback which gives the
object $\opname{Pairs}$ is the fibered product $\opname{Pairs} = Mor
\times_{Obj} Mor$, which gives the usual interpretation as a set of
composable pairs.  This object exists in general since
$\catname{Th(Cat)}$ contains all such finite limits.  If $F$ preserves
finite limits, $F(\opname{Pairs})$ will be a pullback again.
\end{remark}

In Section \ref{sec:wkmod} we briefly recall how a (small) double
category is a model of $\catname{Th(Cat)}$ the \textit{theory of
  categories} in $\Cat$, which is the category whose objects are
(small) categories and whose morphisms are functors.  However, since
we are really interested in \textit{weak} cubical $n$-categories, we
will use the finite limits theory describing \textit{bicategories},
and thus encodes this notion of weakness.

The theory of bicategories, $\catname{Th(\Bicat)}$, is more
complicated than $\catname{Th(Cat)}$, but having seen
$\catname{Th(Cat)}$ we can abbreviate its description somewhat.

\begin{definition} The \textit{theory of bicategories}
$\catname{Th(\Bicat)}$ is the category with finite limits generated by
  the following data:
\begin{itemize}
\item Objects: $\opname{Ob}$,  $\opname{Mor}$,  $\opname{2Mor}$
\item Morphisms:\\
  Source and target maps:
  \begin{itemize}
  \item $s,t:\opname{Mor} \rightarrow \opname{Ob}$
  \item $s,t:\opname{2Mor} \rightarrow \opname{Ob}$
  \item $s,t:\opname{2Mor} \rightarrow \opname{Mor}$
  \end{itemize}
  Composition maps:
  \begin{itemize}
  \item $\circ : \opname{MPairs} \rightarrow \opname{Mor}$
  \item $\circ : \opname{HPairs} \rightarrow \opname{2Mor}$
  \item $\cdot : \opname{VPairs} \rightarrow \opname{2Mor}$
  \end{itemize}
  where
  \begin{itemize}
  \item $\opname{MPairs} = \opname{Mor} \times_{\opname{Ob}} \opname{Mor}$
  \item $\opname{HPairs} = \opname{2Mor} \times_{\opname{Mor}} \opname{2Mor}$
  \item $\opname{VPairs} = \opname{2Mor} \times_{\opname{Ob}} \opname{2Mor}$ 
  \end{itemize}
  are given by pullbacks as in $\catname{Th(Cat)}$
%  \begin{equation}
%    \xymatrix{
%        &  & \opname{Mor} \ar[dr]^{t} & \\
%      \catname{MPairs} \ar[r]^{i} & \opname{Mor}^2 \ar[ru]^{\pi_1}
%        \ar[rd]_{\pi_2} & & \opname{Ob} \\
%        &  & \opname{Mor} \ar[ur]^{s} & \\
%    }
%  \end{equation}
%      (and similarly for $\opname{HPairs}$ and $\opname{VPairs}$)
  \begin{itemize}
  \item The \textbf{associator} map
  \[
    a : \opname{Triples} \rightarrow \opname{2Mor}
  \]
 such that $a$ makes the following diagram commute:
  \[
    \xymatrix{
     \opname{Pairs} \ar[d]^{\circ} &  \opname{Triples} \ar[r]^{1 \times \circ} \ar[d]^{a} \ar[l]_{\circ \times 1}  & \opname{Pairs} \ar[d]^{\circ} \\
     \opname{Mor} & \opname{2Mor} \ar[r]_{t} \ar[l]^{s} & \opname{Mor} \\
    }
  \] (and additional diagrams with the interpretation that $a$ gives
 \textit{invertible} 2-morphisms).
  \item \textbf{unitors}
  \[
    l,r : \opname{Ob} \rightarrow \opname{Mor}
  \] with the obvious conditions on source and target maps.
  \end{itemize} 
\end{itemize}
This data is subject to the usual conditions, including composition
rules for 2-morphisms similar to those for morphisms in
$\catname{Th(Cat)}$, as well as the fact that the compositions for
2-morphisms satisfy the interchange law, associator is subject to the
Pentagon identity, and the unitors obey certain unitor laws.
\end{definition}

The preceding being terse, we note that the Pentagon identity for a
bicategory (i.e. a bicategory in $\catname{Sets}$) is generally
described by saying that the two ways for a composite of associators
from $f \circ (g \circ (h \circ k)))$ to $(((f \circ g) \circ h) \circ
k)$ are equal.  We can express this condition formally, in any
category with pullbacks, building from composable quadruples.  The
pentagon identity may be expressed by a commuting diagram which is
given in \cite{dblbicat}, though that paper omits explicit mention of
$\catname{Th(Bicat)}$.  There are similar diagrams for unitor laws.
Commutativity of these diagrams is imposed in $\catname{Th(Bicat)}$.

We do not propose here to explicitly discuss the theories of
tricategories (though see Gordon, Power and Street \cite{gps}, and the
appendix in Gurski \cite{gurski}), tetracategories (though see Trimble
\cite{tt}), and so forth.  Explicit descriptions of these theories
become quite elaborate very quickly.  These and various other higher
categorical structures have been discussed extensively elsewhere by
Cheng and Lauda \cite{chenglauda}, and Leinster \cite{leinster,
  leinster2}.  We do note, however, that, at least the usual, fairly
well understood, definitions of globular $n$-categories these
definitions can be cast as finite limits theories
$\catname{Th(Tricat)}$ and $\catname{Th(Tetracat)}$, and so on.  So
the process we shall describe here can be applied to all these
theories, giving structures which are correspondingly weaker.  With
bicategories, we have the first case where ``weak'' is meaningful.

Having described the finite limits theory of bicategories, we consider
how to use these theories to define cubical $n$-categories.

\section{Models As Cubical $n$-Categories}\label{sec:models}

We have described the theories of categories and bicategories.  Our
idea here is to see how cubical $n$-categories can be built by taking
models of them in the appropriate setting.

\subsection{Models of Categories and Bicategories}\label{sec:wkmod}

\begin{definition}
A \textit{model} of a finite limits theory $\catname{T}$ in a category
$\catname{C}$ with finite limits is a finite limit-preserving functor
\[
F : \catname{T} \rightarrow \C
\]
\end{definition}

To begin with, we will consider $\catname{T} = \catname{Th(Cat)}$, and
see how to describe the operation of taking ``$n$-fold categories'', a
form of strict $n$-category.

A (small) category is a model of the theory $\catname{Th(Cat)}$ in
$\catname{Set}$.  That is, it is a functor $F : \catname{Th(Cat)}
\rightarrow \catname{Set}$, which is specified by choosing a set of
objects and a set of morphisms, together with set maps making these
into a category.  In this setting, the pullback construction means
that when the target of a morphism $f$ is the source of $g$, there is
a composite morphism $g \circ f$ from the source of $f$ to the target
of $g$ just as expected.

The theory of categories encodes the usual category axioms in terms of
commuting diagrams in $\catname{Th(Cat)}$.  Thus the axioms, in the
presentation of a theory, amount to imposing relations between the
arrows.  For instance, the axiom for associativity can be expressed by
the commuting diagram:
\begin{equation}
  \xymatrix{
    \opname{Triples} \ar[d]_{\circ} \ar[r]^{\circ \times id} & \opname{Pairs} \ar[d]^{id \times \circ} \\
    \opname{Pairs} \ar[r]_{\circ} & Mor
  }
\end{equation}
Since $\catname{Th(Cat)}$ has all finite limits, there are objects
denoting composable $k$-tuples of morphisms for each $k$, similar to
$\opname{Pairs}$, such as
\begin{equation}
  \opname{Triples} = \opname{Mor} \times_{\opname{Ob}} \opname{Mor} \times_{\opname{Ob}} \opname{Mor}
\end{equation}

This is a model of $\catname{Th(Cat)}$ in $\Set$.  Our motivating idea
is to consider models of $\catname{Th(Cat)}$ in the target category
$\C = \Cat$.  Such a model $F$ gives a category $\catname{Ob} =
F(Obj)$ of ``objects'' and a category $\catname{Mor} = F(Mor)$ of
``morphisms'', with functors $s$ and $t$, $Id$, and $\circ$ satisfying
the usual category axioms.  Note that these axioms give conditions at
both the object and morphism level, in addition to those which follow
from the fact that they are functors. Functoriality means that there
are compatibility conditions between the categorical structures in the
two directions.  In fact, these amounts to precisely the definition of
a \textit{double category}.  The ``horizontal'' category is
$\catname{Ob}$ and the ``vertical'' category consists of the objects
in $\catname{Ob}$ and $\catname{Mor}$ together with the object maps
from the functors $s$, $t$, and so forth.  The square 2-cells of the
double category are the morphisms of $\catname{Mor}$.  It can readily
be checked that this gives the usual notion of a double category.
(This is discussed in Leinster \cite{leinster}).

Moreover, a natural transformation $\nu : F \ra G$ between two models
$F, G :\catname{Th(Cat)} \ra \Cat$ is just a double functor in the
usual sense.  In particular, there are functors $\nu(Ob) : F(Ob) \ra
G(Ob)$ and $\nu(Mor) : F(Mor) \ra G(Mor)$, and similarly for $Pairs =
Mor \times_{Ob} Mor$ and so on.  The object and morphism maps of each
of these give assignments for the objects, horizontal and vertical
morphisms, and squares of the double category.  The fact that $\nu$ is
natural gives compatibility conditions between all these maps and the
relevant $s$, $t$, $Id$ and $\circ$ which give exactly the fact that
these assignments define a double functor (in particular, that there
is a functor between the vertical categories).

This describes a \textit{strict} model of $\catname{Th(Cat)}$ in
$\Cat$.  A \textit{weak} model would satisfy the equations in the
category axioms only up to a 2-morphism in $\Cat$, namely up to
natural transformation.  As before, there are categories
$\catname{Ob}$ and $\catname{Mor}$, and functors $s$, $t$, $Id$, and
$\circ$.  However, the equations which hold for categories are
identities in $\catname{Th(Cat)}$, such as associativity, which are
mapped to 2-morphisms - that is, natural transformations in $\Cat$.
That is, regarding the category $\catname{Th(Cat)}$ as having identity
2-morphisms, a weak model allows equations (identity 2-morphisms) to
map to non-identity natural transformations in $\Cat$.  To ensure
coherence, this must be done so that any diagrams of such
2-isomorphisms must commute.  MacLane's coherence theorem (see
e.g. \cite{maclane}, for the result in the context of monoidal
categories) implies that it is sufficient to have the pentagon
identity and unitor identity to imply commutativity of all such
diagrams.

So, for instance, composable pairs would be defined by \textit{weak}
pullback (in $\Cat$) rather than strict pullback (as in $\Set$), so
that in the diagram (\ref{eq:pairs}), instead of satisfying $t \cdot
\pi_1 \cdot i = s \circ \pi_2 \cdot i$, there would only be a natural
isomorphism $\alpha: t \cdot \pi_1 \cdot i \rightarrow s \circ \pi_2
\cdot i$.  Such a weak model is the most general kind of model
available in $\Cat$, but this does not give a general weak cubical
$n$-category.  In particular, it composition is weak in only one
direction. This is equivalent to the definition of a
\textit{pseudocategory} (see, for instance, \cite{pseudocategory}, and
again \cite{leinster}).

In particular, the reason we have weaker composition rules in one
direction for a pseudocategory than the other, from this point of
view, is that taking the target category $\C = \Cat$ determines that
the horizontal structures are (strict) categories, while weakening the
axioms from $\catname{Th(Cat)}$ implies we have vertical bicategories
since, for instance, equations for associativity are mapped to
associator isomorphisms which satisfy the pentagon identity.

Pseudocategories have a number of natural examples when one has two
different types of morphism between the same objects, and in one case,
composition is naturally defined only up to isomorphism.  For example,
there is a pseudocategory whose objects are sets, and where the
horizontal category has functions as morphisms, and the vertical
category has \textit{spans} of sets, composed by pullback.  A related
example has rings as objects, homomorphisms as horizontal morphisms,
and bimodules (composed by tensor product) as vertical.  Square cells,
in these examples, consist of maps of the spans, or bimodules, that
are compatible with the horizontal maps.

Unfortunately, pseudocategories are only weak in one direction, and
strict in the other.  Moreover, for reasons of well-formedness, it is
impossible to use squares as the 2-morphisms to weaken composition in
both directions.  Yet in general, we would like a definition which is
weak in both directions---in the conclusion we return to a class
of examples where this is the natural choice.  We return
to this in the conclusion.  For now, we look at this definition.

\subsection{Internal Bicategories}\label{sec:internal}

The fact that a pseudocategory contains ``vertical bicategories''
suggests a generalization of our approach to double categories.  This
is to consider (strict!) models of $\catname{Th(\Bicat)}$ in
$\Bicat$---that is, functors $F : \catname{Th(\Bicat)} \rightarrow
\Bicat$.  It will be relatively straightforward to treat
$F(\opname{Obj})$ as a horizontal bicategory, and the objects of
$F(\opname{Obj})$, $F(\opname{Mor})$ and $F(\opname{2Mor})$ as forming
a vertical bicategory, but we note that a diagrammatic representation
of, for instance, 2-morphisms in $F(\opname{2Mor})$ would require a
4-dimensional diagram element.  These structures, termed
\textit{double bicategories}, are described in \cite{dblbicat}.

Though we have defined them in terms of the theory of bicategories, as
with double categories, these structures can also be described in
elementary terms.  They have nine types of components, namely the
objects, morphisms, and 2-morphisms in each of the bicategories
$F(\opname{Obj})$, $F(\opname{Mor})$ and $F(\opname{2Mor})$.  There
are a number of connecting ``face maps'' derived from the $s$ and $t$
maps, composition rules, and so on.  The most natural way to express
these diagrammatically involves cells of dimension up to 4, drawn as
products of 0-, 1-, and 2-cells.  For example, in both horizontal and
vertical directions, there are 3-dimensional morphisms like the
``pillow'' $P$ here:
\begin{equation}\label{xy:pillow}
  \xymatrix{
    x \ar[r] \ar[d]^{f'}="1" \ar@/_2pc/[d]_{f}="0" & y \ar[d]^{g'} \\
    x' \ar[r] \uriicell{F_1} & y' \\
    \ar@{=>}"0" ;"1"^{\alpha}
  } \qquad \Rightarrow_P \qquad
  \xymatrix{
      x \ar[r] \ar[d]_{f} & y \ar[d]_{g}="3" \ar@/^2pc/[d]^{g'}="2"  \\
      x' \ar[r] \uriicell{F_2} & y' \\
      \ar@{=>}"3" ;"2"^{\beta}
  }
\end{equation}
Diagrammatically, $P$ should be drawn as the product of an edge and a
(globular) 2-cell.  The vertical pillows are the morphisms of
$F(\opname{2Mor})$, while the horizontal pillows are the 2-morphisms
of $F(\opname{Mor})$

Since one naturally might hope for a fully weak cubical 2-category to
have cells of dimension at most 2, it is also convenient that certain
double bicategories, satisfying ``niche-filling'' conditions, give
rise to what Dominic Verity previously called double bicategories, and
we call ``Verity double bicategories''.  These have horizontal and
vertical bicategories, as well as squares like a double category,
whose composition laws in both directions are weakly associative (up
to 2-cells, rather than squares as in pseudocategories).  We will
return to this in Section \ref{sec:niche}.  For now, if we instead
take the definition of a weak cubical 2-category by internalization as
natural, and follow it, we can see how to extend it to $n$-categories.

Now, a model of a theory $\catname{T}$ as a functor $F : \catname{T}
\ra \mathcal{C}$, hence the maps between models are natural
transformations $\nu : F \ra G$.  In particular, $\nu$ gives, for each
object $t \in \catname{T}$, a morphism $\nu(t) : F(t) \ra G(t)$; for
every morphism in $\catname{T}$, there is a naturality square.  So, in
particular, this defines a concept of morphism of models, and thus the
category of all models, which we denote $[\catname{T} \ra
  \mathcal{C}]$.  In the case of a simple algebraic theory such as
$\catname{Th(Grp)}$, this defines the notion of group homomorphism in
any given setting (say, continuous homomorphism between topological
groups, if $\mathcal{C} = \catname{Top}$.

Now again we consider the strict case.  For $\catname{T} =
\catname{Th(Cat)}$, such a $\nu$ defines, in particular, maps
$\nu(\opname{Ob}) : F(\opname{Ob}) \ra G(\opname{Ob})$ and
$\nu(\opname{Mor}) : F(\opname{Mor}) \ra G(\opname{Mor})$, which
commute with the maps of $\catname{T}$ (in the sense of commuting
naturality squares).  As we saw, such a $\nu$ defines the notion of a
functor between categories, internal to $\mathcal{C}$.  In particular,
if $\mathcal{C} = \Cat$, this defines the usual notion of a functor
between double categories, in the language of models in $\Cat$.

So consider the functor category $\mathcal{C} = [\catname{Th(Cat)},
  \Cat]$, whose objects are double categories (models of
$\catname{Th(Cat)}$ in $\Cat$), and whose morphisms are double
functors (natural transformations $\nu$).  The functor category
inherits the property of having all finite limits from the target
$\Cat$, since the limit of a diagram of functors, at each object $X
\in \catname{Th(Cat)}$, gives the limit of the diagram applied to $X$.
So we can take this new $\mathcal{C}$ as our new target category for
models of our finite limits theory.  Then models $F :
\catname{Th(Cat)} \rightarrow \mathcal{C}$ are triple categories: that
is, cubical 3-categories.  Natural transformations between such models
are triple-functors in a natural sense, and this gives a new category.
Iterating this process gives the usual notion of strict cubical
$n$-category as an $n$-fold category.

Ihe analogous process for $\catname{T} = \catname{Th(\Bicat)}$ and
$\mathcal{C} = \Bicat$, the (1-)category whose objects are
bicategories and whose morphisms are homomorphisms (2-functors)
between them.  Then double bicategories may be seen as functors $F :
\catname{T} \ra \mathcal{C}$, and natural transformations between
these models give a notion of double bifunctor.  Then there is a
functor category $[\catname{T}, \mathcal{C}]$.  Taking this to be our
new $\mathcal{C}$, and iterating the process, we get a notion of weak
cubical $n$-category for all $n$ as an $n$-tuple bicategory.  In
particular, we can inductively define:

\begin{definition}
A weak cubical $0$-category is just a set, and a functor between these
is a set function, so we say $0-\catname{tupleBicat} = \Set$.  Given a
category $(n-1)-\catname{tupleBicat}$, define $n-\catname{tupleBicat} =
  [\catname{Th(Bicat)}, (n-1)-\catname{tupleBicat}$.
\end{definition}

To describe these more completely, note that we can inductively
describe all the types of data which make up an $n$-tuple bicategory,
given by the model $F$.  In particular, there will be a natural
interpretation where we have $3^n$ types of morphism (including
objects as 0-morphisms).  This is beacuse there is one type of data
(elements) when $n=0$; and for $n>0$, we have $(n-1)$-tuple
bicategories $F(\opname{Ob})$, $F(\opname{Mor})$ and
$F(\opname{2Mor})$, each with $3^{n-1}$ types of data.
Diagrammatically, these can be represented as $n$-fold products of
dot, edge, and globular 2-cell (indicating which type we select at
each step of the induction).  Thus, the morphisms naturally have
dimension up to $2n$, although they are composable in only $n$
directions.  Each dimension corresponds to one stage of the inductive
construction.

The composition rules follow from the fact that the composition in any
direction is given by the bicategory axioms in the corresponding stage
of the construction.  For instance, composition in each direction for
cubes (products of edges) is weakly associative in each direction,
where the associator is a morphism given as a product of $(n-1)$ edges
(in all other directions), and a 2-cell (in the direction of
composition).  On the other hand, composition of 2-cells in a
bicategory (i.e. within the category $hom(x,y)$) is strict.  So
composition of those morphisms given as products of a 2-cell with
other data will be strict in the corresponding direction.  However, it
makes sense to say that an $n$-tuple bicategory is a ``weak cubical
$n$-category'', since at least the cubes have weakly asssociative
composition in all $n$ directions.

We have remarked that the term ``double bicategory'' was used by Verity to describe a somewhat different structure, in which all morphisms are naturally represented as 2-dimensional.  Since it seems reasonable to expect that an $n$-category should have morphisms of dimension at most $n$, we briefly consider how the two are related.

\subsection{Niche-filling Conditions}\label{sec:niche}

It is an impetus for much research that there are various relations
between different definitions of $n$-category.  In particular,
relations between (strict) cubical and globular $n$-categories have
been described, by Brown and Higgins \cite{brownhiggins}, among
others.  So for example, double categories are related to
2-categories: if horizontal and vertical morphisms can be composed,
then squares can be considered to be 2-cells between the composites of
the edges.  For higher $n$, and weaker notions of $n$-category, more
complex relations become possible.  In \cite{dblbicat}, we discussed
``niche-filling conditions'' which reduced a double bicategory in the
sense of a model of $\catname{Th(Bicat)}$ in $\Bicat$, to a Verity
double bicategory.  Such a structure has horizontal and vertical
bicategories, as well as squares, together with various composition
rules, and also actions of 2-morphisms on squares.  Here we briefly
give an account of these niche-filling conditions in terms of models
of a finite-limits theory.

A niche-filling condition states that, given any suitable combination
of cells (that is, objects, morphisms, 2-morphisms, and so on) of a
particular shape, there is some other cell which ``fills the niche''
by completing the diagram in a specified way.  A simple example of a
niche-filling condition is fulfilled by the composition operation for
morphisms in a category (i.e. model of $\catname{Th(Cat)}$ in $\Set$).
Here, the niche is given by a choice of object $x$ and pair of
``composable'' morphisms $f$ and $g$ with $s(g) = x$ and $t(f) = x$.
Then the filler for this niche is the composite $g \circ f : s(f) \ra
t(g)$.

Such conditions play a major role in definitions of simplicial
$n$-categories (generalized to Joyal's ``quasicategories''
\cite{joyal-quasi}, also called $\infty$-categories by Lurie
\cite{lurie}) in which the various axioms for a category amount to
``horn-filling conditions''.  For example, composition is replaced by
a 2-simplex (a 2-morphism in a simplicial $n$- or $\infty$-category),
or rather the condition that, for morphisms $f$ and $g$ with $t(f) =
s(g)$, there is a triangle $C$ and edge $g \circ f$ filling the
diagram:
\begin{equation}
\xymatrix{
\drriicell{F} & y \ar[dr]^{g} & \\
x \ar[ur]^{f} \ar@{-->}[rr]_{g \circ f} &   & z
}
\end{equation}
Other properties for such $\infty$-categories are also expressed as
horn-filling conditions

In the cubical case, the more complex shapes of the niches and the
greater number of distinct operations make the situation slightly
trickier, but they are also less crucial to our chosen definition.
They do, however, give modifications to it.  An example of a useful
niche-filling condition in a double category $\mathcal{D}$ would be
the following.  Suppose that, given a horizontal arrow $f$ and
vertical arrow $g$ in $\mathcal{D}$, where the source of $g$ is the
target of $f$, there is a unique invertible square $F$ and vertical
arrow $h$ making the following commute:
\begin{equation}
\xymatrix{
 X \ar[r]^{f} \ar[d]_{h} \driicell{F}  & Y \ar[d]^{g} \\
 Z \ar[r]_{\opname{id}} & Z
}
\end{equation}
Then one can define $h$ to be the composite $g \circ f$, and get a
category generated by the morphisms of both $\Hor$ and $\Ver$, where
the $F$ in the above is regarded as the identity, and all other
squares are discarded.  In fact, to do this, we do not necessarily
need that there is a unique such ``niche-filler'', only that there is
a specified way to find one, which we can use to define composition
(and that these choices are coherent).

There are analogous conditions for double bicategories: horizontal and
vertical ``action conditions'', and a compatibility condition, which
turn the intrinsically four-dimensional structure into a 2-dimensional
structure satisfying Verity's definition of ``double bicategory''.  As
with double categories in the example above, a bicategory can be
obtained from a double bicategory by allowing composition of
horizontal and vertical morphisms.  Here, we are less interested with
this, than with imposing conditions which give a new notion of
$n$-category from an old one, by discarding certain higher morphisms
(in this case, the 3- and 4-dimensional ones), in a consistent way.
In particular, the condition of interest specifies \textit{actions of
  2-cells on squares}.  Given a 2-cell $\alpha$ and a square $F_1$
which share an edge (i.e. are composable), the condition allows one to
complete the left half of the ``pillow'' diagram (\ref{xy:pillow})
with three data.  These are a (unique, invertible) 3-dimensional cell
$F$, the opposite square $F_2$, taken as the ``composite'' $\alpha
\circ F_1$, and with the 2-cell $\beta = \opname{id}$.

In general, a niche-filling condition demands, given a certain
collection of cells (i.e. data from the final model in the chain),
that there are other cells which compose with them, satisfying some
commutation conditions.  A strong version of a niche-filling condition
requires that such niche-filling data are unique.  A weaker version
merely requires that some filler exist.  The sort of condition we want
is one which specifies a filler given a niche: that is, we require
that our weak $n$-category be equipped with maps specifying the
fillers of any niche.  This includes the case where there is a unique
filler.  Given such fillers, we can define actions of one type of
morphism on another by taking the fillers between them to be ``thin''.
That is, we consider them to be the identity, and discard all other
morphisms of the chosen shape.

Given such a niche-filling condition, we can consider the category of
all models which are equipped with such a map.  Provided these
categories have finite limits (i.e. that the niche-filling condition
is preserved by taking finite limits), we can use them in our process
of taking iterated models.  Starting with
$[\catname{Th(Bicat)},\Bicat]$, the category of double bicategories in
our sense, we can find a category consisting of models in this
category \textit{equipped with} a map giving niche-fillers of the kind
used in the three action conditions.  Each of these determines a
double bicategory in the sense of Verity.  Functors which preserve the
niche-filling maps make these into a category $\catname{VDB}$.  In
fact, $\catname{VDB}$ contains finite limits, so we can take
$[\catname{Th(Bicat)},\catname{VDB}]$.

Now, a model $F \in [\catname{Th(Bicat)}, \catname{VDB}]$ determines Verity
double bicategories $F(Ob)$, $F(Mor)$, and $F(2Mor)$ and the
connecting double functors such as $F(\circ)$, and so forth.  This can
be interpreted as a (weak) cubical 3-category.  Moreover, it is weaker
in the new direction, since there are higher-dimensional cells
representing, for example, squares in $F(2Mor)$, which would be
four-dimensional.  If we want our notion of weak $n$-category to have
morphisms represented by cells of dimension at most $n$, we again need
a niche-filling conditions here which would specify which of these
cells to use when pasting.

That is, we would need to specify cells of various dimensions which
define the niche to be filled.  One approach, then, is at each step to
take the full subcategory of models which satisfy these conditions.
However, the number of conditions grows at each step, since there are
more directions in which composites need to be defined.  A different
approach would be to incorporate the niche-filling conditions into our
theory.  However, this means inverting the process used so far, in
which our target category $\mathcal{C}$ for the model $F$ changes, but
the theory $\catname{Th(Bicat)}$ stays the same.  To address
niche-filling conditions at the level of the theory, we would need to
obtain $\catname{Th(DblBicat)}$, a theory of double bicategories
(namely, of functors from $\catname{Th(Bicat)}$ into $\Bicat$), and
add specified maps which give the niche-fillers, satisfying the
implied conditions.

Neither of these approaches to reducing the dimension of the cells of
our weak $n$-categories is particularly elegant, so here we will adopt
the view that the fully general definition is in some sense simpler.

\section{Conclusion}

The presentation we have given for weak cubical $n$-categories
suggests the simple definition that they are $n$-fold bicategories:
models given by applying the functor $[ \catname{Th(Bicat)}, - ]^n$ to
$\Set$.  Thus, $[\catname{Th(Bicat)}, \Set] = \Bicat$ by definition.
Then $[\catname{Th(Bicat)},\Bicat]$ can be denoted
$\catname{DblBicat}$, $[\catname{Th(Bicat)},\catname{DblBicat}]$ can
be denoted $\catname{TrplBicat}$ and so forth.  But note that at each
step of the iteration, we might make a difference choice of which
theory to model in the category produced at the previous step, and of
how strict the model should be.  In particular, we have not considered
the question of what a weak model of $\catname{Th(Bicat)}$ in $\Bicat$
would be, but treating $\Bicat$ as a mere category ignores its full
structure.  $\Bicat$ is most naturally a tricategory in which the
morphisms are 2-functors, 2-morphisms are natural transformations, and
3-morphisms are modifications, allowing weak models.  So in fact, this
schema can generate many different definitions of cubical $n$-category
of many different degrees of strength and weakness.  We are
restricting attention to strict models, since the structures these
produce are already weak enough for some relevant applications.

It is not unusual for different applications to suggest different
definitions of $n$-category, which is one reason for their abundance
(see \cite{chenglauda}).  In particular, $n$-fold bicategories are
quite natural for extending the classes of examples discussed in
\cite{dblbicat} based on (co)spans (also developed extensively by
Grandis \cite{grandis1, grandis3}).  This is a generalization of the
bicategory of cospans in a category $\mathcal{C}$ with pushouts.
Given $X,Y \in \mathcal{C}$, the morphisms in
$\catname{Cospan(\mathcal{C})}$ from $X$ to $Y$ are diagrams $X
\rightarrow S \leftarrow Y$, and 2-morphisms between cospans are
``cospan maps'' given by morphisms $f : S \rightarrow S'$ in
$\mathcal{C}$ which commute with the maps from $X$ and $Y$.  Spans
compose by pushout along two common inclusions.  In the cubical case,
one treats $n$-fold products of such diagrams.  This naturally fits
the framework of the $n$-fold bicategories described here.

A special case is the topological example discussed in
\cite{dblbicat}, which was the main motivation there.  Here, the
cospans of interest are \textit{cobordisms} between manifolds.  That
is, $S$ is a manifold with boundary, and the maps from $X$ and $Y$ are
inclusions of boundary components (in the smooth case, the extra
structure of a collar is needed, which is discussed in a general
setting by Grandis).  Among other things, cobordisms give a way to
study manifolds (and in particular find invariants for them) by
factoring them into pieces, dealing with each piece, and then
composing the pieces.  Typically, the category $\catname{nCob}$ is
described as having $(n-1)$-dimensional manifolds as objects, and
\textit{diffeomorphism classes} of cobordisms as morphisms.  A
diffeomorphism of a cobordism, fixing the boundary (and its collar, if
it has one), is just a cospan map in $\catname{Cosp(Man)}$.  Now, in
particular, taking cobordisms, not equivalence classes, as morphisms
means that composition (by gluing cobordisms at boundary components)
is only weakly associative.

This is why a weak structure was desired.  The need for a cubical
$n$-category comes from a generalization motivated by applications to
(topological) quantum field theories.  This is a generalization to
cobordisms \textit{between} cobordisms---in particular, cobordism
between manifolds with boundary (suitable for field theories in
backgrounds with boundary conditions).  In particular, as a double
bicategory, the objects are $(n-2)$-dimensional manifolds, the
(horizontal and vertical) morphisms are $(n-1)$-dimensional
cobordisms, and the squares are $n$-dimensional cobordisms with
corners.  However, there are other types of morphism here.  Horizontal
and vertical 2-cells are diffeomorphisms of horizontal and vertical
cobordisms.  There are also diffeomorphisms of the $n$-dimensional
body.  Those which fix (pointwise) the horizontal source are our
``vertical pillows'', and vice versa.  Those which fix only objects
are the top-dimensional cells in the double bicategory.

This can of course be extended to a $k$-fold bicategory of cobordisms
with corners having codimension $k$ between the objects,
$(n-k)$-dimensional manifolds, and the top-level morphisms, the
$n$-dimensional cobordisms with corners, up to and including the case
when $n=k$.  In each case, we again have cobordisms as the
cubical-shaped morphisms, and diffeomorphisms (fixing various
components of boundaries) as the remaining morphisms.  Each type of
morphism appearing in our $n$-fold bicategories has a natural
intepretation.

It is possible, even convenient, to discard some of the complexity of
the double bicategory by taking $n$-cobordisms only up to
diffeomorphism (this guarantees the niche-filling conditions discussed
earlier, and yields a Verity double bicategory, which itself can be
further reduced to a bicategory).  But the most natural way to
describe the full structure of cobordism is to leave in these
morphisms.  This means the most natural organizing structure is the
one perhaps most simply understood through the process of
internalization we have described here---that is, in terms of models
of finite limits theories.

\end{document}